\documentclass [12pt]{article}
\input{artadj.sty}
\usepackage{latexsym}
\usepackage{graphicx}

\usepackage{amsfonts}
\usepackage{amsmath, amssymb}

\renewcommand{\thesection}{\arabic{section}}

\newcommand{\limn}{\lim_{n\rightarrow\infty}}

\newcommand{\limtau}{\lim_{\tau\rightarrow\infty}}

\newcommand{\goto}{\rightarrow}

\newcommand{\bo}{{\bf 1}}

\newcommand{\beqa}{\begin{eqnarray*}}
\newcommand{\eeqa}{\end{eqnarray*}}

\newcommand {\beq} {\begin{equation}}
\newcommand {\eeq} {\end{equation}}
\newcommand {\bear} {\begin{eqnarray}}
\newcommand {\eear} {\end{eqnarray}}
\newcommand {\bears} {\begin{eqnarray*}}
\newcommand {\eears} {\end{eqnarray*}}   

\newcommand {\done} {\quad\vrule height4pt WIDTH4PT}

\newcommand {\barr} {\begin{array}}
\newcommand {\earr} {\end{array}}

\newtheorem{thm}{Theorem}[section]
\newtheorem{lem}[thm]{Lemma}

\newtheorem{cor}[thm]{Corollary}

\graphicspath{%
    {converted_graphics/}
    {/}
}
\begin{document}
\bibliographystyle{plain}

\title{On {\em BASTA} for Discrete-Time Queues}
\date{}

\author{ Muhammad El-Taha \\
                  Department of Mathematics and Statistics\\
                       University of Southern Maine\\      
                      96 Falmouth Street\\
                            Portland, ME  04104-9300\\
                     Email:el-taha@maine.edu
}
%
\maketitle

\ls{1}

\noindent{\bf Abstract.} We  address  the  issue  of   Bernoulli/Geometric  arrivals  see  time-averages, {\em BASTA}, in  discrete-time  queues  when  using  a  variety  of  scheduling  rules.  It  is  well-known  that  {\em BASTA}/{\em ASTA} holds  in  the  framework  of  a  discrete-time  process  with  an imbedded  arrival  process  without  regard  to  scheduling  rules.  However,  when scheduling  rules are invoked,  it  is  not  always  clear  whether  arrivals  see  time  averages  as  asserted  by  the  general  theory.  We  resolve  this  issue  so  that  {\em BASTA}/{\em ASTA}  holds  regardless  of  the  scheduling  rule  at  play.
\bigskip

\noindent
{\bf Keywords:} Discrete-time queues, {\em BASTA}, {\em ASTA}, {\em GASTA}
\bigskip

\noindent

\newpage
\ls{1.30}
\section{Introduction}\label{sec:intro}

In discrete-time systems, {\em BASTA} (Bernoulli arrivals see time averages)  is a fundamental property that states that, when arrivals follow a Bernoulli process the stationary distribution as seen by arrivals  equals  the time-average stationary distribution. 
In discrete-time queues, time is divided into slots of equal lengths. We assume w.l.o.g. that slots are  one-unit length, where  slot boundaries (edges) are numbered by $\tau=1,2,\ldots$. Because  arrivals and departures are assumed to occur at slot boundaries one needs to identify the order of arrivals and departures and if the state of the system is observed before or after an arrival and/or a departure. The physical behavior of the actual system  leads to various scheduling rules like {\em EAS} (early arrival system) and {\em LAS} (late arrival system).
 We give a brief description of the five most cited rules in the literature in Subsection 2.1.
Clearly, there is a gap between  the general theory, of {\em BASTA}, based on a discrete-time process with an imbedded arrival process (when the imbedded process is Bernoulli) and the application of this theory to discrete-time queues with various scheduling rules. In continuous time  queueing systems, a time-average distribution function refers to a distribution function obtained by averaging continuously over time. In  discrete-time queues, a time-average distribution function refers to averaging over the set of slots edges that  will be referred to
as the random observer epochs. Contrary to continuous time queues, one can also utilize the distribution function  at slot centers, which is referred to 
 as the outside observer epochs. Additionally, in discrete-time queues, there are other points of reference that one can utilize.  
{\em This article addresses this issue by providing a simple resolution that  asserts that the time-average distribution function should be averaged at  {\bf \em potential} pre-arrival epochs. The rest of the article is a formalization of this notion.}

The literature  on discrete-time queues  contains  several instances of the gap between the theory of {\em ASTA} (arrivals see time averages) as a fundamental property  and it application to discrete-time queues.
Consider  a discrete-time process with an imbedded discrete point (arrival) process and the {\em LAA} (lack of anticipation  assumption) or its equivalent.
On one hand, there are results in the literature (\cite{Hal83},~\cite{Mak89},~\cite{Elt99}, and~\cite{Miy92}) that assert that   if the arrival point process is Bernoulli then the stationary distribution function at  pre-arrival instants and the time-average stationary distribution function coincide. On the other hand, there are examples of discrete-time queues (\cite{Gra92},~\cite{Dad01}, and~\cite{Des02})  where arrivals are Bernoulli, but pre-arrival probabilities do not equal the time-average probabilities in the sense  that time averaging takes place at slot edges. 
Halfin~\cite{Hal83} is the first to  address {\em BASTA}, which he refers to as {\em GASTA} (geometric arrivals see time averages) in the framework of a  discrete process  with an imbedded  Bernoulli arrival  process. His result and proof mimics the statement and proof of the corresponding continuous time {\em PASTA} (Poisson arrivals see time averages) given by ~\cite{Wol82}. Later~\cite{Mel90,Mel90b},~\cite{Sti89}, and~\cite{Elt92,Elt99,Elt02} among others generalize {\em PASTA} to {\em ASTA} by removing the requirement that arrivals follow a Poisson process and weakening the {\em LAA}  assumption into the {\em LBA} (lack of bias assumption). Moreover,~\cite{Mak89} discuss {\em ASTA} for discrete time processes. They show that in discrete-time, {\em ASTA} holds under the weaker discrete time {\em LBA} condition.

In applications to discrete-time queues the state of the system (number of customers) can be observed at a number of special epochs, most prominently are the  random observer epochs or slot edges  and the outside observer epochs or slot centers. 
Using {\em BASTA} as given by~\cite{Hal83}, Gravey and Hebuterne~\cite{Gra92} show that the pre-arrival probability distribution is identical to  the outside observer's distribution function  for the {\em LAS}  rule and point out that the same result does not hold for the {\em EAS}  rule. In other words   averaging over the set of slot centers, {\em BASTA}  holds for the  {\em LAS}  rule but not for the {\em EAS} rule.
Daduna~\cite{Dad01} addresses the {\em BASTA} issue in Chapter 2 and states,``in discrete time, however, a {\em PASTA} analogue usually does not hold". He also concludes that averaging at slot edges;  the pre-arrival stationary distribution  does not always equal the stationary equilibrium  distribution  even if the arrival process is Bernoulli. Desert and Daduna~\cite{Des02} study instances when {\em BASTA} holds and conclude that the order of  events (arrivals and departures) play a crucial rule in obtaining {\em BASTA} results. 
 A sample of recent articles that use {\em BASTA}/{GASTA} include  \cite{Bar19,Bar20,Cha22,Wan21,Cha19}, among others.


This article is organized as follows:
In Section~\ref{sec:ps} we give a brief introduction of the five scheduling rules  used  in discrete-time queueing systems. To motivate our approach we give  brief review of {\em BASTA} focusing on Halfin's approach.  
In Section~\ref{sec:pasta} we  focus on a characterization of  {\em ASTA}  and {\em BASTA} as it applies to discrete-time queues, and prove the main results and discuss consequences. We also give  give applications for the late-arrivals departures-first ({\em LA-DA}) rule as it presents the most challenging situation. Additionally, we give  {\em BASTA} results for birth-death queues.

\section{Background and Motivation}\label{sec:ps}

Several authors have address discrete-time {\em BASTA} in the framework of a discrete-time process with an imbedded process, notably~\cite{Hal83,Mak89,Elt92,Elt99,Elt02}.  
 In this section, we  briefly explain  discrete-time scheduling rules, present the system setup, and give  a brief commentary    to point out the gap when {\em BASTA} is applied to discrete-time queues.

\subsection{Discrete-Time Queues with Scheduling Rules}
 We give a brief description of the five most cited rules in the literature.
Depending on the behavior of the actual system, the order of potential arrivals and departures at any given slot vary significantly.

\noindent
{\bf {\em EAS} Model.}
In the early arrival system ({\em EAS}), in a given slot,  potential arrivals are scheduled to occur before potential departures.
More specifically a potential arrival in  slot $(\tau,\tau+1]$  occurs in $(\tau,\tau+)$, and a potential departure in slot $(\tau -1,\tau]$ occurs  in $(\tau-,\tau)$.  Moreover, if an arrival finds an idle server, it goes into service immediately  and can potentially depart in the same time slot.


\noindent
{\bf {\em LAS-IA} Model.}
In the late arrivals system ({\em LAS}) the order of potential arrivals and departures is reversed so that  potential departures occurs early in a time slot and potential arrivals occur at the end of the slot.  More specifically, a potential departure in  slot $(\tau,\tau+1]$  occurs in $(\tau,\tau+)$, and a potential arrival in slot $(\tau -1,\tau]$ occurs  in $(\tau-,\tau)$.
Moreover, if an arrival finds an idle server and goes into service immediately,  and can potentially depart at the start of the following  time slot, the system is called immediate access ({\em LAS-IA}).

\noindent
{\bf {\em LAS-DA} Model.}
 This scheduling  rule is similar to the  {\em LAS} rule  except that when an arrival finds an idle server it waits until the start of the next slot to start service, thus the system is called delayed access ({\em LAS-DA}).



\noindent
{\bf {\em LA-AF} Model.} In this late arrivals with arrivals-first rule, both potential arrivals and departures occur at the end of a slot with potential arrivals taking place before potential departures.

\noindent
{\bf {\em LA-DF} Model.} In this rule,  potential arrivals and departures occur at the end of a slot, but now  with potential departures taking place before potential arrivals.

These rules are depicted in Figure 1. For additional information about different scheduling regimes, one may consult~\cite{Hun83b} ~\cite{Cha00},~\cite{Gra92},~\cite{Elt97a} and~\cite{Dad01}.

\begin{figure}[tbp] 
  \centering
  \includegraphics[bb=0 0 1376 769,width=5.67in,height=3.17in,keepaspectratio]{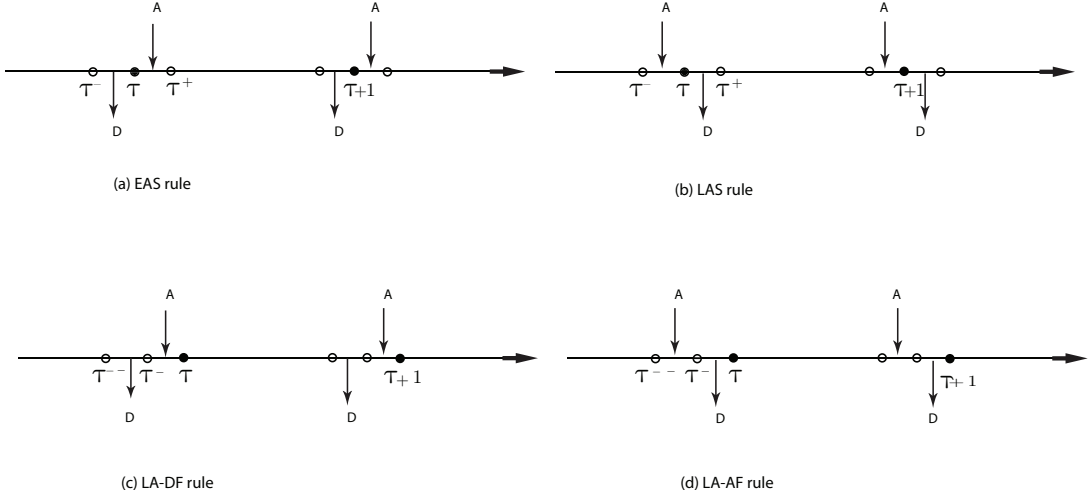}
  \caption{Scheduling rules for discrete-time queues}
  \label{fig:SRfig}
\end{figure}

\subsection{System Setup}

To formalize the system setup, let $\{Z(\tau), \tau\geq 1 \}$ be a 
discrete-time process with countable state space $S$.  We assume that $S$ is a complete separable metric space with the  Borel $\sigma$-field $\cal{B}(S)$.
Since we
shall be working  with sample-path properties, it is
helpful to think of $Z=\{Z(\tau), \tau\geq 1 \}$ as a deterministic
process.  In the context of a stochastic model, $\{Z(\tau), \tau\geq 1 \}$ should 
be interpreted as a particular sample path
(realization) of the corresponding stochastic process.
The process makes a transition from one state to another, possibly itself, at every time instant. 
In other words, we allow $Z(\tau)$ to make a
transition from a state to itself. 
Additionally, let  $\{N(\tau), \tau=1,2,\ldots \}$  be a discrete-time point process such that $N({\tau})$ is the number of points (arrivals) at time $\tau$. 
This definition allows for  multiple points, i.e. more than one arrival at the same time,(e.g.  batch arrivals). For a simple point process $N(\tau)\in\{0,1\}$. 

\subsection {Classical {\em BASTA}.}

Let $\{N(\tau),\tau=1,2,\ldots \}$ be a discrete Bernoulli process where $N(\tau)=0,1$ with $0$ indicates no arrival and $1$ indicates an arrival at time $\tau$.
Let $B\in\cal{B}(S)$ be any set in the state space of $Z$.
As in~\cite{Hal83}  define the following limits when they exist
\begin{eqnarray*}
\pi(B)&=&\limtau \sum_{k=1}^{\tau}\bo\{Z(k)\in B\}/\tau\;;\\
\pi^A(B)&=&\limtau \sum_{k=1}^{\tau}N(k)\bo\{Z(k)\in B\}/\sum_{k=1}^{\tau}N(k)\;.
\end{eqnarray*}
In a queueing system $\pi(B)$ and $\pi^A(B)$ represent, respectively, the time-average and pre-arrival average probability distributions.
Moreover, the {\em LAA} states that: 
 For every $n$, the set of random variables  $Z(1), Z(2),\cdots,Z(n)$ is independent of the $N(n), N(n+1),\cdots$.
 
Halfin~\cite{Hal83} proves that under {\em LAA}, for any set of states $B$, $\pi^A(B)=\pi(B)$.
Note that with these definitions it is {\em implicit} that the time-average probabilities $\pi(.)$ are averaged at the slot boundaries $\tau$, and potential arrivals occur right after the slot boundary (this is necessary to insure $N(n)$ and $Z(n)$ are independent as stipulated by the {\em LAA}).

 Potential arrivals and departures are typically scheduled to take place either right before or right after slot boundaries $\tau$. 
Here, the state of the system is observed at the slot boundary, so one needs to specify whether $N(\tau)$ takes place before or after $Z(\tau)$ is observed. The set-up of a discrete time process with an imbedded Bernoulli process implicitly fits with the {\em EAS} rule but not with the {\em LAS} rules. Moreover, one can re-define the reference points to occur at slot centers and with little care can insure that the {\em LAA} hold thus in this situation the arrival distribution function coincides with the outside observer distribution function as done by~\cite{Gra92}.

Note that in queueing systems where $Z(\tau)$ represents the number of customers in the system, {\em LAA}  assumes that $N(\tau)$ and $Z(\tau)$ are independent. However,
 $Z(\tau)$ depends on arrivals up to $\tau$,  so that it is implicitly assumed that an arrival at $\tau$ should be scheduled after $\tau$, i.e. in $(\tau,\tau+)$. 
This would rule out a scheduling rule like {\em LAS}.
In~\cite{Mak89} they define 
$
\pi^A(B)=\limtau \sum_{k=1}^{\tau}N(\tau)Z(\tau-1).
$
This definition insures that an arrival at $\tau$ is observing the system at $\tau-1$.  This works fine in systems
 where it is implicitly assumed that no departures occur in $(\tau-1,\tau)$, but
in queues with scheduling rules  where a potential
  departure takes place after $\tau-1$ and before $\tau$,  what an arrival sees  will depend on whether an actual departure occurred in that time slot. This is the case for {\em EAS} and {\em LA-DA} systems.
Similar issues arise when applying the results of~\cite{Mak89}, and~\cite{Elt99} to discrete-time queues with certain scheduling rules. 

Interestingly  depending on the scheduling rule, Bernoulli  arrivals see time averages at slot edges, slot centers, or other reference points. 
In Section~\ref{sec:pasta} we  present our main results that  resolve the issues presented in this section.

\section{{\em ASTA} with Bernoulli Arrivals}\label{sec:pasta}

In this section, we show that in discrete-time queues with any scheduling rule, Bernoulli  arrivals always see time averages  obtained at {\em potential} arrival epochs. This is remarkable as it bridges the gap between existing  theory on  {\em BASTA} and it applications to discrete-time queues with scheduling rules.

\subsection{Characterization of {\em ASTA}}

Let $\tau^A$ be a {\em potential} arrival epoch. Note that for each discrete time point $\tau$, there is a $\tau^A$ associated with it and this association depends on the scheduling rule.  For $\tau=1,2,\ldots$,  $\tau^A \in (\tau,\tau^+)$ for {\em EAS} models, $\tau^A \in (\tau^-,\tau)$  in {\em LAS}  and {\em LA-DF} models, and  $\tau^A \in (\tau^{--},\tau^-)$ for {\em LA-AF} models. 
Moreover, let  $\{\tau^A_n, n=1,2,\ldots\}$ be an associated deterministic point process,  not necessarily simple. That is  there can be more than one arrival. We shall refer to $\tau^A_n$ as an arrival event.
Additionally, let  $A(\tau)=\max\{n; \tau^A_n\leq \tau\}$ be the number of arrival events in $(0,\tau]$. In the case of a simple  arrival process   $A(\tau)$ counts the number of arrivals in $(0,\tau]$.
Furthermore,   $\{N(\tau^A), \tau=1,2,\ldots \}$  is a discrete-time point process such that $N({\tau^A})$ is the number of points (arrivals) at time $\tau^A$. 
This definition allows for  multiple points, i.e. more than one arrival at the same time, (e.g.,  batch arrivals). For a simple point process $A(\tau^A)= \sum_{u=1}^{\tau^A} N(u^A)$. 

Let $p(N(\tau^A)=k)=p(k), k=0,1,\ldots$, $\sum_{k=0}^{\infty}p(k)=1$, with finite mean $\lambda'=\sum_{k=0}^{\infty}kp(k)$. 
Here, $Z(\tau)$ can be thought of as number of customers in the system at time $\tau$.
With this notation, the {\em LAA} assumption will take the following form
\medskip

\noindent
{\bf Lack of Anticipation Assumption ({\em LAA}).}
Assume that 
\begin{equation}\label{eq:LAA2}
P(N(\tau^A)=k|Z(m^A-)=z_m; 1\leq m^A\leq\tau^A)=p(k); 
\end{equation}
for all $\tau^A\geq 1$, $z_1,\ldots, z_{\tau^A}; k=0,1,\ldots$, and $\{p(k)\}$ is a p.m.f.

 For any $B\in \cal{B(S)}$ define the following limits when they exist
\begin{eqnarray*}
\lambda &:=& \limtau A(\tau)/\tau =\limn n/\tau^A_n;\\
\lambda (B) &:= &\limtau \sum_{k=1}^{A(\tau^A)} \bo\{Z(\tau^A_k-)\in B\}/\sum_{u=1}^{\tau} \bo\{Z(u^A)\in B \}\;;\\
\pi(B)&:=&\limtau\sum_{u=1}^{\tau} \bo\{Z(u)\in B\}/\tau\;;\\
\pi^{A}(B)&:=&\limn \sum_{k=1}^n \bo\{Z(\tau^A_k-)\in B\}/n \;;\\
\pi^{PA}(B)&:=&\limtau \sum_{u=1}^{\tau} \bo\{Z(u^A)\in B\}/\tau^A \;;\\
\pi^O(B)&:=&\limtau \sum_{u=1}^{\tau} \bo\{Z(u-.5)\in B\}/\tau\;.
\end{eqnarray*}
We interpret $\lambda $ as the the asymptotic arrival event  frequency, $\lambda (B)$ as the  asymptotic  state $B$-arrival frequency, $\pi(B)$  as the long-run frequency of time in  state $B$, $\pi^A(B)$ as the long-run frequency of arrivals that find the system in set $B$, $\pi^{PA}(B)$  as the long-run frequency of time in  set $B$ where the system is observed at potential arrival instants, and $\pi^{O}(B)$  as the long-run frequency of time in  set $B$ where the system is observed at slot centers. Note that if the point process is simple then  w.p.1, $\lambda $ and $\lambda (B)$ are, respectively,  the state-independent and state-dependent arrival probabilities.

Note that $\pi^{PA}(.)$ are similar to $\pi(.)$ except we replace slot boundaries $\tau$  by the corresponding potential arrival epochs.  These epochs take place right before or right after  slot edges $\tau$ depending on the scheduling rule. Contrary to random observer epochs, the potential arrival epochs are scheduling rule dependent. 
Now, we give the following result that is the discrete-time analogue of the covariance  formula given by ~\cite{Elt99}. 

\begin{thm}\label{thm:cov} Consider a discrete-time queueing system with any scheduling rule. For any set of states $B\in \cal{B}(S)$  assume that limits $\lambda$ and $\lambda(B)$ exist, then  $\pi(B)$ exists if and only if  $\pi^{PA}(B)$ exists, in which case
\begin{equation}\label{eq:cov}
\lambda (B)\pi^{PA}(B)=\lambda \pi^{A}(B)\;.
\end{equation}
\end{thm}
\noindent
{\bf Proof.}
We use the discrete-time version  of  $Y=\lambda X$, see El-Taha and Stidham~\cite{Elt99}. Let $Y(\tau^A)= \sum_{k=1}^{A(\tau^A)} \bo\{Z(\tau^A_k-)\in B\}$, that is $Y(\tau^A)$ is the number of points (arrival events) that occur in the interval $(0,\tau^A]$. Now, $X_k=\sum^{k}_{n=1}\bo\{Z(\tau^A_n-)\in B\}- \sum^{k-1}_{n=1}\bo\{Z(\tau^A_n-)\in B\}=\bo\{Z(\tau^A_k-)\in B\}$, so that
 $\sum_{k=1}^m X_k= \sum^{m}_{k=1}\bo\{Z(\tau^A_k-)\in B\}$, and  $\pi^{A}(B)= X$.
  Using  $Y=\lambda X$, 
$\pi^{A}(B)= Y/\lambda$, so that
\begin{equation}\label{eq:1}
\lambda \pi^{A}(B)=\limtau \frac{\sum_{k=1}^{A(\tau^A)} \bo\{Z(\tau^A_k-)\in B\}}{\tau^A}\;.
\end{equation}
Now, let $Y(B,\tau^A)=\sum_{u=1}^{\tau} \bo\{Z(u^A)\in B \}$, then
\begin{eqnarray}
\lambda(B)\pi^{PA}(B) &=&  \limtau \frac{\sum_{k=1}^{A(\tau^A)} \bo\{Z(\tau^A_k-)\in B\}}{Y(B,\tau^A)}\frac{Y(B,\tau^A)}{\tau^A} \nonumber\\
   &=&
\limtau \frac{\sum_{k=1}^{A(\tau^A)} \bo\{Z(\tau^A_k-)\in B\}}{\tau^A}\;. \label{eq:2}
\end{eqnarray}
The proof follows by comparing (\ref{eq:1}) and (\ref{eq:2}).
\done

Theorem~\ref{thm:cov}, serves as a starting point in proving $ASTA$, and {\em BASTA}.
 For the next result we assume that  all relevant limits are well-defined and that  $\pi(B)$ and $\pi^{A}(B); B\in \cal{B}(S)$ are proper frequency distributions. 

\begin{cor}[Conditions for {\em ASTA} in Discrete-Time Queues]\label{cor:astadisc}
 Suppose $0 <\lambda(B) < \infty$.
Then, for each $B\in \cal{B}(S)$,

\begin{equation}\label{eq:asta}
 \pi^A(B) = \pi^{PA}(B) \mbox{ if and only if }
 \lambda(B) = \lambda.
\end{equation}

Moreover, the following are equivalent:

(i) $\lambda(B)$ is independent of $B \in \cal{B}(S)$;

(ii) $\lambda(B)= \lambda$  for all $B \in  \cal{B}(S)$;

 (iii) the frequency distributions $\pi^A(B); B\in \cal{B}(S)$ and $\pi^{PA}(B); B\in \cal{B}(S)$  coincide.
\end{cor}

\noindent
{\bf Proof.}
 The proof follows from Theorem~\ref{thm:cov}.
\done
\medskip

The equivalence of (ii) and (iii) says
that the time-average frequency distribution at  pre-arrival instants  coincides with the
frequency distribution just before 
potential arrivals   if and only if the conditional arrival frequency in
each set of states $B$ coincides with the unconditional  frequency. In a stochastic framework, sometimes it is easier to work with assumptions weaker than the {\em LAA} to prove this equivalence.  The lack of dependence assumption
({\em LDA}) states that the events $(N(\tau^A)=1)$ and $(Z(\tau^A-)=n)$ are independent. The lack of bias assumption ({\em LBA}) states that the two events are uncorrelated.

The {\em LAA} implies that the {\em LBA} holds, i.e. $\lambda (B)=\lambda$, so that by Corollary~\ref{cor:astadisc}
\[
\pi^A(B)= \pi^{PA}(B) \;\; w.p. 1\;;  
\]
for every $B\in\cal{B}(S)$.
See ~\cite{Mak89} and~\cite{Elt99} for details.

Note that these results are valid for any discrete system with Bernoulli arrivals where the {\em LAA} assumption holds. This includes single server as well as multi-server models. 
The following consequences follow immediately from this theorem.
\begin{cor}[{\em BASTA}]\label{cor:basta1}
Consider a queueing system having Bernoulli arrivals with parameter $\alpha$, $0<\alpha <1$.  Assume that arrivals are independent of  system state, i.e. the  {\em LAA} holds. Then
for all $n\geq 0$, $P(N(\tau^A)=1|Z(\tau^A-)=n)=\alpha$, $0<\alpha<1$, and for any scheduling rule 
\[
\pi^A(n)=\pi^{PA}(n)\;.
\]
\end{cor} 

This result states that Bernoulli arrivals always see time-averages in the sense that the stationary distribution at pre-arrival epochs will equal the time-average stationary distribution at {\em potential} arrival epochs regardless of the scheduling rule. It remains to characterize the stationary distribution, when possible,  at potential arrival instants for each scheduling rule. Our next result does exactly that.

\begin{cor}[{\em BASTA}]\label{cor:basta2} Under the conditions of Corollary~\ref{cor:basta1}, 
we have  
 $$P(N(\tau^A)=1|Z(\tau^A-)=n)=\alpha,\mbox{\;\;for all\;\;} n\geq 0;. \mbox{ Moreover,}
$$

(i) for the  EAS  models 
\[
\pi^A(n)=\pi(n)\;; n\geq 0\;,
\] 

(ii) for the {\em LA-AF} models 
\[
\pi^A(n)=\pi(n)=\pi^O(n)\;;n\geq 0\;, \mbox{and}
\]

(iii) for the {\em LAS-IA} and {\em LAS-DA} models
\[
\pi^A(n)=\pi^O(n)\;;n\geq 0\;.
\]
\end{cor}
{\bf Proof.} For the  {\em EAS} and {\em LA-AF} models, the state of $Z$ at  epoch $\tau^A-$ coincides with the state at  $\tau$, so  that the distribution functions at these epochs coincide. Moreover, for the  {\em LA-AF} model, the state of $Z$ at  epoch $\tau^A-$ coincides with the state at  $\tau-.5$. This proves (i) and (ii).
For {\em LAS} models, note that  the system state $Z(.)$  does not change in the interval $((\tau-1)^+, \tau-)$, so  the state of $Z$ at $\tau^A-$ and $\tau-.5$ coincide. This shows  that $\pi^{PA}(B)=\pi^{O}(B)$.\done
 
Corollary~\ref{cor:basta1} and Corollary~\ref{cor:basta2} clarify an important {\em BASTA} issue. 
 In Corollary~\ref{cor:basta1}, we give one unified {\em BASTA} result for all scheduling models. Moreover, in Corollary~\ref{cor:basta2} we show how {\em BASTA} leads to arrivals seeing time averages at different reference points  when {\em EAS}, {\em LAS} and {\em LA-AF} scheduling is involved. More specifically, Corollary~\ref{cor:basta2} states that the   pre-arrival distribution function for all scheduling rules ({\em EAS}, {\em LAS-DA}, {\em LAS-IA}, and {\em LA-AF}) is either equal to the random or outside observer's  distribution function except for the   {\em LA-DF} rule, which we address  next.

Now we deal with the {\em LA-DF} rule for  a one-dimensional generalized birth-death discrete-time model where arrival and service completion probabilities, $\alpha(j)$ and $\beta(j)$ are state dependent and state $j$, $j\geq 0$, is observed at pre-potential arrival epochs. 
 Because departures occur before arrivals the number of customers an arrival sees depends on whether an actual departure occurs in $(\tau^{--},\tau^{-})$.  This leads to the following result.
\begin{lem} For any LA-DF generalized birth-death model with state $n$ arrival and service  probabilities $\alpha (n)$ and $\beta (n)$ respectively, we have
\begin{equation}\label{eq:las11}
\pi^A(n)= \frac{1}{\alpha}\left[\alpha(n)(1-\beta(n))\pi(n) +\alpha (n+1)\beta(n+1)\pi(n+1)\right]\;;
\end{equation}
where  $\alpha=\sum_{k=0}^{\infty} \lambda(k)\pi(k).$\\
In particular, if $\alpha(n)=\alpha$ for all $n \geq 0$ (state independent), then
\begin{equation}\label{eq:las12}
\pi^A(n)= (1-\beta(n))\pi(n) +\beta(n+1)\pi(n+1)\;.
\end{equation}
\end{lem}
{\bf Proof.}
If the state at a potential  pre-arrival instant $\tau-$ is $n$, i.e. an arrival sees $n$ customers in the system, then the state at  $\tau$ is $n+1$. So
\begin{eqnarray*}
\pi^A(n)&=&\limtau p(Z(\tau)=n+1|N(\tau-) =1)\\
        &=&\limtau \sum_{r \in S}p(Z(\tau)=n+1,Z(\tau-1)=r,N(\tau-) =1)/p(N(\tau-)=1)\\
 &=&\limtau \sum_{r \in\{n,n+1\}}p(Z(\tau)=n+1|Z(\tau-1)=r,N(\tau-) =1)\\
&& \times p(N(\tau-)=1|Z(\tau-1)=r) p(Z(\tau-1)=r)/p(N(\tau-)=1)\\
       &=& \left[\alpha(n)(1-\beta(n))\pi(n) +\alpha (n+1)\beta(n+1)\pi(n+1)\right]/
\sum_{k=0}^{\infty} \alpha(k)\pi(k)\;,
\end{eqnarray*}
where in the last step, as $\tau \goto \infty$, we used 
$p(N(\tau-)=1|Z(\tau-1)=r)=\alpha(r), r=n,n+1$;
 $p(Z(\tau-1)=r)=\pi(r), r=n,n+1$;
$p(Z(\tau)=n+1|Z(\tau-1)=n,N(\tau-) =1)=1-\beta(n)$;and 
$p(Z(\tau)=n+1|Z(\tau-1)=n+1,N(\tau-) =1)=\beta(n)$.
The lemma follows by noting that $\alpha=\sum_{k=0}^{\infty} \alpha(k)\pi(k)\;$.\done
\medskip

\noindent
{\bf Remark.} Daduna \cite{Dad01} uses this type of argument. Using (\ref{eq:las11}), we see that {\em ASTA} does not hold here in the sense that pre-arrival probabilities equal time-average probabilities where the average is taken with respect to  time instants $\tau$.
In this approach $\beta(n)$ can be evaluated when service times are geometric with parameter $\beta$ in which case $\beta (n)=\beta, n \geq 0$ which, using (\ref{eq:las12}), leads to
\[
\pi^A(n)=(1-\gamma)\gamma^n\;; n=0,1,\ldots \;,
\] 
where $\gamma=\frac{\alpha(1-\beta)}{\beta(1-\alpha)}$. 
Next, we address {\em BASTA} for generalized birth-death processes.

\subsection{{\em BASTA} for Birth-Death Discrete-Time Queues}

Consider a generalized birth-death process 
 with state dependent arrival and service completion probabilities, $\alpha(j)$ and $\beta(j)$ where $j\geq 0$, is observed at pre-potential arrival epochs. Now assume Bernoulli arrivals, i.e. $\alpha(j)=\alpha$, $j\geq 0$, and $LAA$, so that {\em BASTA} holds. Depending on the scheduling rule $\beta(0)$ will take one of two values:  $\beta(0)=0$ or  $\beta(0)=\beta$, $ 0<\beta<1$. For each $j\geq 0$, let
\[
\gamma(j)=\frac{\alpha(1-\beta(j))}{\beta(j+1)(1-\alpha)}\;;
\]
then it follows 
 that for all $n\geq 0$ we have
\begin{equation}\label{eq:bdg}
\pi^{PA}(n)=\left\{\begin{array}{lll}
        \gamma(0)\prod_{j=1}^{n-1}\gamma(j)\pi^{PA}(0),& n \geq 1\;,\\\\
        \left[1+\gamma(0)(1+\sum_{k=2}^{\infty}\prod_{j=1}^{k-1}\gamma(j))\right]^{-1},& n=0\;,
\end{array}  \right.
\end{equation}
where the $\prod(.)$ over an empty set is one.
 This leads to the following result.
\begin{thm}\label{thm:BD}
Consider a birth-death  queue with Bernoulli arrivals with arrival probability $\alpha$, LAA, and state dependent  service completion probabilities $\beta(j);j\geq 0$. Then

(i) For LAS-DA and LA-AF, $\beta(0)=0$ 
and
\begin{equation}\label{eq:bdb=0}
\pi^{A}(n)=\pi^{PA}(n)=\left\{\begin{array}{lll}
        \frac{\alpha}{\beta(1)(1-\alpha)}\prod_{j=1}^{n-1}\gamma(j)\pi^{PA}(0),& n \geq 1\;,\\\\
        \left[1+\frac{\alpha}{\beta(1)(1-\alpha)}(1+\sum_{k=2}^{\infty}\prod_{j=1}^{k-1}\gamma(j))\right]^{-1},& n=0\;.
\end{array}  \right.
\end{equation}

(ii) For EAS, LAS-IA, and LA-DF, $\beta(0)=\beta$ 
and
\begin{equation}\label{eq:bdb=b}
\pi^{A}(n)=\pi^{PA}(n)=\left\{\begin{array}{lll}
        \frac{\alpha(1-\beta)}{\beta(1)(1-\alpha)}\prod_{j=1}^{n-1}\gamma(j)\pi^{PA}(0),& n \geq 1\;,\\\\
        \left[1+\frac{\alpha(1-\beta)}{\beta(1)(1-\alpha)}(1+\sum_{k=2}^{\infty}\prod_{j=1}^{k-1}\gamma(j))\right]^{-1},& n=0\;.
\end{array}  \right.
\end{equation}
\end{thm}
{\bf Proof.} Note that  $\beta(0)=0$ implies that $\gamma(0)=\alpha/\beta(1)(1-\alpha)$, so (i) follows by substituting in (\ref{eq:bdg}) and simplifying. The proof of (ii) is similar.\done.

Theorem~\ref{thm:BD} applies to discrete-time queues with Bernoulli arrivals and general service times. In this case however, 
$\beta(j)$ is not readily available except when service times are geometric in which case $\beta(j)=\beta, j\geq 1$. Moreover, 
using Theorem~\ref{thm:BD}, one can say {\em LAS-DA} and {\em LA-AF} arrivals see one time-average distribution while {\em EAS}, {\em LAS-IA}, and {\em LA-DF} see a different time-average distribution.


%

\noindent
{\bf Acknowledgments.}
The author thanks the editors and reviewers for comments that led to improvements in the presentation of this article.
This research did not receive any specific grant from funding agencies in the public, commercial, or not-for-profit sectors.


\begin{thebibliography}{10}

\bibitem{Bar19}
F.~P. Barbhuiya and U.C. Gupta.
\newblock Discrete-time queue with batch renewal input and random serving
  capacity rule: {GI}$^{X}$/{Geo$^Y$/1}.
\newblock {\em Queueing Systems}, 91:347--365, 2019.

\bibitem{Bar20}
F.P. Barbhuiya and U.C. Gupta.
\newblock A discrete-time {GI}$^{X}$/{Geo}/1 queue with multiple working
  vacations under late and early arrival system.
\newblock {\em Methodology and Computing in Applied Probability}, 22:599--624,
  2020.

\bibitem{Cha22}
M.~Chaudhry and V.~Goswami.
\newblock The {G}eo/{G}$^{a,Y}$/1/{N} queue revisited.
\newblock {\em Mathematics}, 10(17):3142, 2022.

\bibitem{Cha00}
M.~L. Chaudhry.
\newblock On numerical computations of some discrete-time queues.
\newblock In W.K. Grassmann, editor, {\em Computational Probability}, pages
  365--408. Springer, Boca Raton, 2000.

\bibitem{Cha19}
M.~L. Chaudhry, J.J. Kim, and A.D. Banik.
\newblock Analytically simple and computationally efficient results for the
  {GI}$^{X}/{G}eo/c$ queues.
\newblock {\em Journal of Probability and Statistics}, 2019, 2019.

\bibitem{Dad01}
H.~Daduna.
\newblock {\em Queueing Networks with Discrete Time scale}.
\newblock Springer-Verlag, Berlin-New York, 2001.

\bibitem{Des02}
B.~Desert and H.~Daduna.
\newblock Discrete time tandem networks of queues: Effects of different
  regulation schemes for simultaneous events.
\newblock {\em Performance Evaluation}, 47(2-3):73--104, 2002.

\bibitem{Elt92}
M.~El-Taha and S.~Stidham~Jr.
\newblock A filtered {{\em ASTA}} property.
\newblock {\em Queueing Systems: Theory and Applications}, 11:211--222, 1992.

\bibitem{Elt99}
M.~El-Taha and S.~Stidham~Jr.
\newblock {\em Sample-Path Analysis of Queueing Systems}.
\newblock Kluwer Academic Publishing, Boston, 1999.

\bibitem{Elt02}
M.~El-Taha and S.~Stidham~Jr.
\newblock Filtration of {{\em ASTA}}: A weak convergence approach.
\newblock {\em Journal of Statistical Planning and Inference}, 100:171--183,
  2002.

\bibitem{Elt97a}
M.~El-Taha, S.~Stidham~Jr., and R.~Anand.
\newblock Sample-path insensitivity of symmetric queues in discrete-time.
\newblock {\em Nonlinear Analysis, Theory, Methods and Applications},
  30:1099--1110, 1997.
\newblock Proc. 2nd World Congress of Nonlinear Analysts.

\bibitem{Gra92}
A.~Gravey and G.~Hebuterne.
\newblock Simultaneity in discrete-time single server queues with {B}ernoulli
  inputs.
\newblock {\em Performance Evaluation}, 14:123--131, 1992.

\bibitem{Hal83}
S.~Halfin.
\newblock Batch delays versus customer delays.
\newblock {\em Bell System Technical Journal}, 62(7):2011--2015, 1983.

\bibitem{Hun83b}
J.J. Hunter.
\newblock {\em Mathematical Techniques of Applied Probability, Volume {II}:
  Discrete-Time Models: Tecniques and Applications}.
\newblock Academic Press, New York, 1983.

\bibitem{Liu19}
R.~Liu, A.~S. Alfa, and M.~Yu.
\newblock Analysis of an {ND}-policy {Geo/G/1} queue and its application to
  wireless sensor networks.
\newblock {\em Oper Res Int J}, 19:449--477, 2019.

\bibitem{Mak89}
A.~Makowski, B.~Melamed and W.~Whitt.
\newblock On averages seen by arrivals in discrete-time.
\newblock {\em IEEE Conference on Decision and Control. Conference Proceedings;
  Tampa, Florida}, 28:1084--1086, Dec. 1989.

\bibitem{Mel90}
B.~Melamed and W.~Whitt.
\newblock On arrivals that see time averages.
\newblock {\em Operations Research}, 38:156--172, 1990.

\bibitem{Mel90b}
B.~Melamed and W.~Whitt.
\newblock On arrivals that see time averages: A martingale approach.
\newblock {\em J. Appl. Probab.}, 27:376--384, 1990.

\bibitem{Miy92}
M.~Miyazawa and Y.~Takahashi.
\newblock Rate conservation principle for discrete-time queues.
\newblock {\em Queueing Systems}, 12:215--229, 1992.

\bibitem{Sti89}
S.~Stidham~Jr. and M.~El-Taha.
\newblock Sample-path analysis of processes with imbedded point processes.
\newblock {\em Queueing Systems}, 5:131--165, 1989.

\bibitem{Wan21}
R.~Wang.
\newblock Queue size distribution on a new {ND} policy {Geo/G/1} queue and its
  computation designs.
\newblock {\em Discrete Dynamics in Nature and Society}, 2021:1--16, 2021.

\bibitem{Wol82}
R.~Wolff.
\newblock Poisson arrivals see time averages.
\newblock {\em Operations Research}, 30:223--231, 1982.

\end{thebibliography}

\end{document}